\input amstex
\documentstyle {amsppt}
\magnification=1200
\parindent 20 pt
\NoBlackBoxes

\define\g{\gamma}
\define\a{\alpha}

\define \lm{\lambda}
\define\n{\nu}

\define\vp{\varphi}
\define\th{\theta}

 \define\CO{\Cal O}
\define \BZ{\Bbb Z}
\define \BC{\Bbb C}
\define \BP{\Bbb P}

\define\tm{\times}

\define\mk{\medskip}

\define\ep{\endproclaim}
\define\edm{\enddemo}

\define \de {\partial }
\define\1{^{-1}}
\define\2{^{-2}}

\define \opi{\overline{\pi}}

\define \oV {\overline {V}}

\define \tS{\tilde{S}}\define \tU{\tilde{U}}\define \tC{\tilde{C}}

\define \tf {\tilde {f}}\define \tx {\tilde {x}}\define \tw {\tilde {w}}

\baselineskip 15pt

\title
{ Cylinders over affine surfaces}
\endtitle
\leftheadtext{T. Bandman and L. Makar-Limanov}
\rightheadtext{Cylinders over affine surfaces}
\author
T. Bandman and L. Makar-Limanov
\endauthor
\thanks The first author  is supported by the 
 Ministry of Science and
Technology, by the Ministry of Absorption, State of Israel,
and the Emmy N\"other Institute for
Mathematics.
The second  author is supported by the NSF grant. 
\endthanks
\subjclass Primary 13B10, 14E09; Secondary 14J26, 
14J50,  14L30,  14D25, 16W50, 
 \endsubjclass
\keywords
 Affine varieties, $\BC-$actions, locally-nilpotent derivations.
\endkeywords

\abstract
For an affine variety $S$ we consider the ring $AK(S),$
which is the  intersection of the rings of constants 
of all locally-nilpotent derivations of the ring $\CO(S).$
We show that 
$AK(S\times\BC^n)=AK(S)$
for a smooth affine surface  $S$
with $H^2(S,\BZ)=\{0\}.$ 
\endabstract
\endtopmatter 
\document
\mk

 \baselineskip 20pt
 \subheading {Introduction}

In this paper we are trying to understand better a ring invariant, which was
introduced in \cite{ML1}, (see also \cite{KML}). 
This invariant was 
used in order to show that non of the Dimca-Koras-Russell threefolds, 
which are related to the linearizing question, is isomorphic
to $\BC^3$ (\cite{ML2}, \cite{KML}, \cite{KKMLR}). 
It also helped to describe the automorphisms
of a surface $x^ny = f(z)$ (see \cite{ML3}) and to give a new proof of the
theorem of S. Abhyankar, P. Eakin, W. Heinzer \cite{AEH}
on cancelations for curves
(in characteristic zero case) and some generalizations of this theorem
(\cite{ML4}).

Let us start with necessary algebraic notions in the generality which
corresponds to our intended setting.

Let $R$ be an algebra over a field $\BC$.
Then a $\BC$-homomorphism $\de$ of $R$ is called a {\it derivation} of $R$
if it satisfies the Leibniz rule: $\de(ab) = \de(a)b + a\de(b)$.
  
Any derivation $\de$ determines two subalgebras of $R$. One is the kernel of
$\de,$ which is usually denoted by $R^{\de}$ and is called the {\it ring of 
$\de$-constants,} by analogy with the ordinary derivative.
The other is $nil(\de)$, the {\it ring of nilpotency of} $\de$.
It is determined by $nil(\de) = \{a \in R|\de^n(a) = 0, n >> 1\}$.
In other words $a \in nil(\de)$ if for a sufficiently large natural
number $n$ we have $\de^n(a) = 0$.

Both $R^{\de}$ and $nil(\de)$ are subalgebras of $R$ because of the
Leibniz rule.

Let us call a derivation {\it locally nilpotent} if $nil(\de) = R$.
Let us denote by $lnd(R)$ the set of all locally nilpotent derivations.

The best examples of locally nilpotent derivations are the partial
derivatives on the rings of polynomials $\BC[x_1, ..., x_n]$. 

The intersection of the rings of constants of all locally
nilpotent derivations of $R$ is the invariant mentioned above. It 
will be called the {\it ring of absolute constants}
and denoted by $AK (R)$.

If $V$ is a complex affine variety and $\CO(V)$ is the ring of regular
functions on $V,$ let us denote $AK(\CO(V))$ by $AK(V)$. A locally nilpotent
derivation of $\CO(V)$ corresponds via exponentiation to a $\BC$ action on $V$
(see \cite{S}). It follows, that if $AK(V)\neq\CO(V),$ variety $V,$
as we see later, is
ruled or ``cylinderlike,'' which means that it contains an open subset
which is a product of affine variety and a complex line $\BC.$
It seems that the invariant $AK(V)$ is especially helpful
 when one tries to compare a variety with $\BC^n.$
E.g. M. Miyanishi (\cite{Mi1}) showed that 
 an affine surface  $V$ with factorial ring $\CO(V)\neq AK(V), $
 is isomorphic to $\BC ^2, $ provided $\CO(V)=\BC.$ So, we think it is   
rather important to learn to compute this invariant.

One of the approaches to this is to find a connection between 
$AK(V \times W)$,
$AK(V)$, and $AK(W)$ where $V$ and $W$ are affine varieties. E.g. 
it is known (see \cite{ML4}), that if $V$ is a curve which is not
an affine line, then  $AK(V \times W) = AK(V) \bigotimes_{\BC} AK(W)$. 
It is also known (see \cite{ML1}) that if $AK(V) = \CO(V),$ 
then  $AK(V \times \BC ^1) = AK(V) \bigotimes
 \BC = AK(V)$.
Any locally nilpotent derivation $\de$ of $\CO(V)$
 or of $\CO(W)$ can be extended 
to a
locally nilpotent derivation of $\CO(V \times W)$
 by $\de(f) = 0$ for any
$f \in \CO(W)$ (of $f \in \CO(V)$ correspondingly). 
So 
it is clear that $AK(V \times W) \subset AK(V) \bigotimes AK(W)$.
In geometric terms it is the following obvious observation. 
If $f \in \CO(V \times W)$ is invariant under all $\BC$-actions
it is also invariant under all $\BC$ actions which are ``lifted"
from the components.

Unfortunately, it is not true in general even when $W = \BC^1,$ that
$AK(V \times W) = AK(V) \bigotimes AK(W),$ which is demonstrated by 
the following example. 

Let surfaces $S_1$ and $S_2$ be defined in $\BC^3$ by equations 
$xy = z^2 - 1$ and $x^2y = z^2 - 1$.
Danielewski (\cite{D}, \cite {K}) showed that the 
cylinders over these surfaces are 
isomorphic and Fieseler (\cite{F}) 
proved that the surfaces themselves are not isomorphic. 

These surfaces were suggested by Danielewski as a
counterexample to the generalized Zariski cancelation conjecture.

In our setting they provide an example of a situation when  
$AK(R) \subsetneq AK(R[x])$. Let $R_1$ and $R_2$ be the rings of 
regular 
functions on $S_1$ and $S_2$ correspondingly. It is easy to find two 
locally nilpotent
derivations on $R_1$ such that the intersection of their kernels 
is just $\BC$. 
Say, take $\de_1$ defined by 
$\de_1(x) = 0,$ $\de_1(y) = 2z,$ $\de_1(z) = x,$
and $\de_2$ defined by
$\de_2(x) = 2z,$ $\de_2(y) = 0,$ $\de_2(z) = y.$

On the other hand it is possible to show that any 
locally
nilpotent derivation of $R_2$ has $x$ in the kernel and that
$AK(R_2) = \BC[x]$ (\cite{ML3}). Thus
 $AK(R_2) = \BC[x] \neq AK(R_2[x]) = AK(R_1[x])
= \BC$.

So it seems rather natural to find conditions on  
a variety which ensure the equality $AK(V) = AK(V \times \BC^n)$.

The goal of this paper is to show that if $S$ is a smooth surface and 
$H^2(S, \BZ) = 0,$ then $AK(S \times \BC ^n) = AK(S).$  

We also would like to state the following

\proclaim{Conjecture} $AK(V) = AK(V \times \BC^n)$ 
if $\CO(V)$ is a factorial ring.\ep

If this, indeed, is true it will advance rather substantially
an understanding of Zarisky cancelation  conjecture 
and the linearizing question for $\BC ^*-$ action on  $\BC ^n.$
\vskip 0.5 cm
\proclaim {Acknowledgements} It is our pleasure to thank Sh.Kaliman
for reading the manuscript and for his very important remarks. \ep
\vskip 0.5 cm
 \subheading {1. Auxiliary Facts}

Assume that a group $G,$ possibly infinite
dimensional, is generated by a finite number of $\BC$-actions
$\{\vp_i\}$, which act algebraically  on a $n$-dimensional
irreducible reduced  affine variety $X$.
This means that for any $i$ there is a regular rational 
map $\vp_i:\BC\times X\to X$ such that

\rom{a)} $\vp_i(z_0,x)=\vp_i^{z_0}(x) $ is an algebraic
 regular automorphism of $X$;

\rom{b)} $\vp_i^{z_0+z_1}(x)=\vp_i^{z_0}\circ \vp_i^{z_1}(x).$

If the group $G$ is algebraic, then, due to the Rozenlicht Theorem
(see, for example, \cite {P-V}), there are two possibilities:
either the a general orbit is Zariski dense in $X$, or 
there exists a $G-$ invariant rational function on $X.$
\footnote {Further on we shall say that a set $W$ 
is dense in the algebraic set $V,$ if $W$ contains a 
Zariski open subset of $V.$}

The following proposition is a generalization of this fact for a 
non-algebraic group.

\proclaim{Proposition 1.1}
If the $G$-orbit $G_{x_0}= \{ y : y=g(x_0),g\in G \} $ of
 a general point $x_0$ is not dense
in $X,$ then there exists a $G$-invariant  
 rational map $\pi: X\to X_G$ of the variety $X$ into  an
irreducible  algebraic variety $X_G$ with $\dim X_G < \dim X$.\endproclaim

\vskip.5cm
We need two Lemmas to prove the Proposition.
\proclaim{Lemma 1.2}
The graph $D_G=\{(x,y)\in X\times X: y=g(x), g\in G\}$ is a dense subset of
 a closed algebraic subset of $X\times X.$
 Moreover, the orbit of a general point $x_0$ is a dense subset 
 of an algebraic subset of $X$.
 \endproclaim

\demo{Proof of Lemma 1.2}
Let word $I=\{i_1\dots i_s\}$ be a word of length $s,$ where $i_k$ 
are
natural numbers. Let us  define a regular algebraic map
$F_I: \BC^s\times X\to X$ by
$$F_I(z,x)=\vp_{i_1}^{z_1}\circ\dots\circ \vp_{i_{s}}^{z_s}(x),$$
where $z=(z_1\dots z_s)\in \BC^s $ and $x\in X.$

For a multiindex $I=\{i_1, ...i_s\}$ 
the graph of this map $\Gamma_I=\{(z,x,y): y=F_I(z,x)\}$
 is a closed subset of $\BC^s\times
X\times X.$ Since this graph is isomorphic to $\BC^s\times X,$
 it is irreducible. 

Denote by $\overline \Gamma_I$ its closure  in
the product  $\Bbb P^s\times\overline X\times\overline X$
of closures $\Bbb P^s,\overline X$ of $\BC^s$ and $X$ respectively.
It is irreducible because  $\Gamma _I$ is irreducible.
Hence its projection
 $\overline Z_I\subset \overline X\times\overline X$ into $\overline
X\times \overline X$ is an
 irreducible closed subset of $\overline X\times\overline X$,
containing the projection $W_I$ of $\Gamma_I$ as a dense subset (since $\Gamma_I$ is
dense in $\overline\Gamma_I).$
Let $Z_I=\overline Z_I\cap(X\times X).$
Then $W_I\subset Z_I\subset\overline Z_I,$\
where $Z_I$ is
 an irreducible closed subset of $X\times
X$ and $W_I$  is its dense subset.
For any two words $I $ and $J$, such that  $J$ contains $I$, \ $W_I\subset W_J$ and $\overline
Z_I\subset \overline Z_J.$
Now from the irreducibility of both $\overline Z_I$ and $\overline Z_J$
 follows  that either
$\overline Z_J=\overline Z_I$ or $\dim \overline Z_j>\dim \overline Z_I.$
Since $\dim Z_I\le 2n,$ it is possible to chose a word $I$ such
 that $k=\dim \overline {Z}_I$ is maximal among all indices $I$.

Consider two cases.

1. Let $k=2n.$ Then $Z_I=\overline X$ and  for a general point $x_0\in X$
 we have $\dim Z_I\cap \{x_0\times X\}=n.$ Thus its projection 
 together with projection of $W_I$ into
second factor $X$ is dense in $X.$  But the projection of $W_I$  
is precisely the orbit 
of a point $x_0$, which shows that in this case the general orbit
 is  dense
 in $X$. 

2. Let $k$ be less then $2n.$
Since $
\overline {Z}_J\subset\overline {Z}_{J\cup I}$
and by the choice of $I,$  
 $\overline {Z}_J\subset \overline{Z}_I $
for any word $J.$
So,  $W_J\subset \overline Z_I,$ for any such word $J$, 
and $D_G =\cup W_J$ is
contained in $\overline Z_I$ as well.
On the other hand, $D_G$ contains $W_I,$ and 
the last  is dense in $\overline Z_I$. Moreover by definition  $D_G
\subset X\times X.$
Thus,  $D_G$ is a dense subset of a closed subset $Z_I$
 of the product $X\times X.$   The orbit of a general point $x_0$
 is a projection of $D_G\cap (x_0\times X)$ into $X$, which is a dense subset
 of the intersection of the  projection of $\overline
 {Z}_I\cap(x_0\times \overline {X})$ into $\overline X$  and $X.$\qed
\enddemo

\proclaim{Lemma 1.3}
If the general  $G$-orbit $G_{x_0}$ is not dense in $X,$ 
then there are rational
functions on the variety $X,$ invariant under the action of the group $G.$
\endproclaim

\demo{Proof of Lemma 1.3}

Let  
 $K$ be an ideal of the functions in, which are equal to zero on the closure
 $\overline D_G\subset
X\times X$ in $X\times X$ of the set $D_G.$
Any such   function $f(x,y)$ has the form
$$f(x,y)=\sum_{i=1}^{n_f} p_i(x)q_i( y),\tag 1$$
where $(x,y)$ are the points of $X\times X$ and  $q_i(y)\in \CO(X), 
\ \ p_i(x)\in \CO(X)$. 
If for all   functions $f(x,y)\in K$ all $p_i(x)=0,$ then
 any pair $(x,y)\in X\times X$ belongs to
$\overline D_G,$ i.e. $D_G$ is dense in $X\times X,$ \
and the general orbit  $D_G\cap\{x_0\times X\}$ is dense in $X$
for a general point $x_0\in X,$  which contradicts the assumptions of 
the Lemma. Thus, $K$ contains non-zero functions.

Let $f_0$ be a  function in $K\setminus
\{0\},$ such that $\n=n_{f_0}=min\{n_f| f\in K\setminus\{0\}\}:$
$$f_0=\sum_{i=1}^{\n} p_i(x)q_i(y).$$

Once more two cases are possible.

1) $\n=1$ and  $f_0(x,y)=p_1(x)q_1(y).$
\flushpar Then $D_G$ is contained in
  $(\{X\times N\})\cup
(\{M \times X\}),$
where $N=\{y|q_1(y)=0\}$
and $M=\{x|p_1(x)=0\}.$ 
For  any point $x_0$
the pair $(x_0,x_0)\in D_G \subset   (\{X\times N\})\cup
\{(M \times X\})$
which is impossible, if $x_0\not\in N\cap M.$
 Thus, $\n\geq 2.$

2) $f_0(x,y)=p_1(x)q_1(y)+\sum_{i=2}^{\n}p_i(x)q_i(y),$
where $p_1(x)\ne 0.$
Then for $(x,y)\in D_G$
 and  $g\in G$ two equalities hold:
 
$$
p_1(x)q_1(y)
+\sum\limits_{i=2}^{\n}p_i(x)q_i(y)
=f_0(x,y)
= 0.$$
(2)
 $$p_1(gx)q_1( y)+\sum\limits_{i=2}^{\n}
p_i(gx)q_i(y)=
f_0(gx,y)=0.$$

Therefore $$p_1(gx)f_0(x,y)-p_1(x)f_0(gx,y)=
\sum\limits_{i=2}^{\n} (p_1(gx)p_i(x)-p_1(x)p_i(gx))q_i(y)\in K$$
and is ``shorter'' then $f_0.$
 Since $\n$ was minimal by the choice of $f_0,$
it means that $p_1(gx)p_i(x)-p_1(x)p_i(gx)\equiv 0$
for $i=2,...,\n.$

Thus
$\dsize\frac
{p_i(x)}{p_1(x)}$ 
are $G-$invariant rational functions. They  cannot be constant: if, say,
$\dsize\frac
{p_2(x)}{p_1(x)}=c,$
 then $f_0$ could have been written in a ``shorter `` way:
$$f_0=p_1(x)(q_1(y)+cq_2(y))+\sum\limits_{i=3}^{\n}p_i(x)q_i(y).$$

Hence   $\dsize\frac
{p_i(x)}{p_1(x)}$  
are the needed G-invariant rational functions.\qed\enddemo

\demo{Proof of the Proposition 1.1}

Let $f_0(x,y)=p_1(x)q_1(y)+\sum\limits_{i=2}^{\n}
p_i(x)q_i( y).$
Then the map $$\pi: x\to(p_1(x):p_2(x):\dots :p_{\n}(x))$$
 is a rational map  of
$\overline X$ into a projective space $\BP^{\n-1}.$
Since there always exists a resolution $X'$
of the map $\pi,$ (i.e. an irreducible projective variety
birationally equivalent to $\overline X,$ such that the induced map
$\pi ':X' \to \pi (\overline X)$
is regular), the set $\pi(\overline X)\subset \BP^{\n-1}$ is an image of 
an irreducible projective set under a regular map.  
Thus, $\pi (\overline X)$ is projective
and irreducible (\cite {Sh2}, 5.2) and contains $\pi (X)$
as a dense subset. Since the general orbit
$G_x$ is dense in a closed subset $\overline G_x$ of $X,$
 $\dim \pi (\overline X)\leq \dim (\overline X) - \dim 
\overline G_x < \dim (\overline X)$, which completes 
the proof. \qed \enddemo

The next Lemma is a 
particular case of Lemma 2.2 in a paper of M. Miyanishi \cite {Mi2}.

  \proclaim {Lemma 1.4} (see \cite {Mi2}) 
  Let $R$ be a finitely generated ring. Then 
  it has a non-zero locally nilpotent derivation 
  if and only if there exists an element $t\in R,$
  such that $R[t^{-1}]$ is isomorphic to a polynomial ring $S[x]$.
\ep

\demo {\it Proof of Lemma 1.4} See \cite {Mi2}.
\qed\enddemo   

   \proclaim {Corollary 1.5} Let $X$ be an affine normal
   variety and $\pi:X\to Y$ be a regular map into
   a normal affine variety $Y.$ 
   Let $t\in \CO(Y)$ and let $D$ be divisor of its zeros.
   Assume that $V=X\setminus\pi^{-1} (D) 
   \cong( Y\setminus D)\times \BC.$ Then there is a $\BC$-action
   on $X, $ such that its general orbit is a fiber of  
   the map $\pi.$ \ep

\demo {\it Proof of Corollary 1.5}
Let $S=\CO(Y\setminus D)$ . Then $\CO(V)=\CO((Y\setminus D)\times \BC)
=S[x].$ On the other hand $\CO(V)=\CO(X)[t^{-1}],$ since for any function
 $r\in \CO(V)$ there is such positive integer $k$ that
$t^kr\in \CO(X).$ So,  
according to Lemma 1.4 there is a locally nilpotent 
derivation $\de$ on $R,$ such that $S\subset R^{\de}.$  
Since variety $Y$ is affine, the functions $s\in S$ 
divide points in $Y$. That means that a general fiber 
of $\pi$ may be described as $s_1=const, .....s_n=const$
for some $s_1,...s_n \in S.$  Since all $s_i$ are $\de-$  constants,
it means that the general fiber is invariant under 
$\BC-$action, corresponding to $\de.$ On the other hand, 
 a general fiber is isomorphic to $\BC$ and consequently
 is an orbit.  
\qed \enddemo.

\subheading {2. Invariant of product}

In this   
section we prove the following 

\vskip .5 cm
\proclaim {Theorem 2.1 } Let $S$ be a smooth surface with  
$H^2(S,\BZ)=\{0\}. $
Then $AK(S\times \BC ^n)= AK(S).$\ep
\vskip .5 cm

\proclaim {Remark} The condition $H^2(S,\BZ)=\{0\} $  is     
essential. In the introduction an example is given of a surface with
$H^2(S,\BZ)\neq\{0\},  $ and  such that     
  $AK(S\times \BC ^n)\neq AK(S).$\endproclaim

Let $S$ be a smooth affine surface, $X=S\times \BC^n $ and $\pi:X\to S$
the natural projection. Assume that there is a $\BC-$action $\vp_{\lm}$
on $X$ such that the  orbit $\Gamma_{x_0}=\{\vp_{\lm}( x_0), 
\lm \in \BC\}$
of a general point $x_0$ is not contained in $\pi\1(\pi(x_0))$.
Denote by $\psi_{1,\lm},...,\psi_{n,\lm}$ the standard actions
acting along the fibers of a projection $\pi$ and by $G$ the 
automorphisms group of $X$ generated by $\vp_{\lm}, \psi_{1,\lm},...
\psi_{n,\lm}$. In Lemma 1.2  we proved that
the orbit $G_{x_0}$ of a general point
$x_0$ is a  dense subset of a 
closed subset $\overline {G_{x_0}}$ of $X$.

 If $\dim \overline {G_{x_0}}=n+1$
then by Proposition 1.1 there exists a dominant 
$G-$invariant rational  map $p:X\to \BP ^1.$

 The fibers of this map contain $G-$ orbits.

Since the map $p$ is $G-$invariant, it induces the 
map $p_1:S\to \BP^1,$
such that the following diagram is commutative

$$\alignat{12}\
& \ &&X && &&\\
&\pi&& \downarrow \searrow p \tag 3 \\
& \ && S \ \ \overset{p_1} \to\longrightarrow\ \BP^1 \endalignat$$

We   consider two different cases.

{\it Case 1}. The map $p_1$ is regular (i.e.
it is everywhere defined).

{\it Case 2}. 
For any 
$G-$invariant rational  map $p$ the map  $p_1: S\to \BP^1$ 
is not regular.

 \proclaim {Lemma 2.2 } In  case 1 there is a 
$\BC -$ action $\th_{\lm}$ on the surface $S,$ such that a general
orbit $\g_{s_0}=\{\th_{\lm}( s_0),\lm \in\BC \}=\pi(G_{x_0})$, where 
$  \pi(x_0)=s_0.$\endproclaim

\demo {\it Proof of  Lemma 2.2}

In this case both maps $p_1$ and $p$ are regular.

Choose a closure $\tS$ of $S$ in such a way that the map
 $p_1$ may be extended to a regular map $\tilde p_1: \tS \to \BP ^1.$

 Let $\tC$ be the 
normalization of the Stein factorization of  a map $\tilde p_1.$
By definition (see, e.g, \cite {B}, p.66)
 that means, that $\tC$ is a smooth 
curve included into the following 
commutative diagram:

$$\alignat{12}\
& \ &&\tS \ \ \overset{\tilde p_1}\to\longrightarrow \BP ^1 \\
& h\ &&\downarrow \nearrow q &&\tag 4 \\
& \ && \tC  && && \endalignat$$
where all the maps are regular, $h$ has connected fibers
 and $q$ is finite.
Let $F_c=h^{-1} (c)\subset \tS$ 
be a  fiber over 
a general point $c\in \tC.$
Take  points 
$s\in F_c\cap S$ and $ x\in \pi\1(s).$ The orbit
$\Gamma_x=\{\vp_{\lm}(x), \lm \in \BC \} \subset X$ is a rational 
curve with a single puncture, because it is an image of a complex
 plane. It 
follows, that $\pi(\Gamma_x) \subset S$ is a rational curve with a single 
puncture as well. Since $F_c$ is a closed connected curve, containing 
$\pi(\Gamma_x),$ $F_c$ has to be the  closure of $\pi(\Gamma_x) .$ 
Moreover, it has to be smooth
( for general $c$), since it is a fiber of a regular map 
of a smooth surface onto a smooth curve 
(see, for example, [Sh1], \S 4). 
Thus, the restriction of $h$
onto $S$ has a general fiber $h\1(c)\cap S= \pi (\Gamma_x),$ 
which is isomorphic to 
the complex plane, and a  general fiber $h\1(c)$ is an irreducible smooth
rational curve. 
Hence,  $\tS$ is a ruled surface, (see [Sh1], Theorem 2, 
chapter 4). Moreover,   the divisor  $D=\tS\setminus S$   
has precisely one irreducible component   $D_0$ which 
is mapped by $h$ isomorphically onto $\tC,$
 since a general fiber of $h_{|S}$
is isomorphic to $\BC^1.$ 
All other components of the divisor $D$
do not intersect with a general fiber of $h,$ which is irreducible. 
Consider the  reducible fibers $F_i, i=1,,,s.$  Since  intersection
$(F_i, D_0)=1$, they have the following structure:
$F_i= C_i+\sum_{j=1}^{n_i}\a_j E_{ij},$
where $C_i, E_{ij}$ are irreducible components, $(C_i,D_0)=1,
 (E_{ij},D_0)=0$. According to
[Mi3], 4.4.1, every reducible fiber contains at least one 
exceptional curve of the first type.

Since $H^2(S)=\{0\},$
the inclusion $D\to \tS$ induces the epimorphism $H_2(D) \to
H_2(\tS).$

The group $H_2(\tS,\BZ)$ may be described as follows.
Let $S_0$ be a surface obtained by   
blowing down  all the exceptional curves  
of the first type
 in all the fibers of the map $h.$ We may repeat this
 procedure till me obtain the       surface $S_n, $  the map $t:\tS\to  S_n,$
and the regular map $h_n: S_n\to\tC$ with irreducible fibers. 
Thus, it will be geometrically ruled surface. 
The group $H_2(S_n, \BZ)$ is a direct sum $d_0.\BZ \oplus
f.\BZ,$ where $d_0 $ and $f$ are the homology classes of $D_0$
 and a general fiber respectively ([B],III.18). Let $E_{ij}$ be all the
 irreducible curves in $\tS,$ which are contracted by $t,$
 and $e_{ij}$
their homology classes ( there is precisely $n_i$ of such curves
in reducible fiber $F_i$).
Then $H_2 (\tS, Z)= H_2(S_n)\oplus \sum_{i} \sum_{j}  e_{ij}.\BZ.$
  Thus, $H_2(D)\to H(\tS)$ 
may be an epimorphism only in case when $D$ contains at least one fiber $F$ 
and other fibers have only one irreducible component in $S.$

\bigskip

 There are two important consequences of this fact. 
\bigskip
1)$p_1(S)$ is an affine subset of $\tC.$   Indeed, a point
$c=h(F)\in (\tC\setminus h(S)).$

2) $S$ has a ``cylinderlike'' subset.
Indeed, since $\tS$ is ruled,  by taking away
the finite number of points $c_1, c_2, ...,c_N$ from $\tC$ we obtain
the Zariski  open subset $\tU \subset \tS,$ which is isomorphic to
$(\tC\setminus\{c_1,...,c_N\}) \tm \BP^1,$ (\cite{B}, p. 26)
and $\tU\cap S$ 
is isomorphic to  $(\tC\setminus\{c_1,...,c_N\}) \tm \BC$ 

 Adding, if needed, some other points, we may assume that
$c_1+...+ c_N$  is the  zero-divisor of a regular function 
on $\tC\setminus h^{-1}(c).$ 

By Corollary 1.5 there is a $\BC-$action $\th_{\lm}(s)$ on
 $S$ such, that an orbit $\g_s=\{\th_{\lm}(s),\lm \in \BC\}$ of a general
point $s$ coincides with $h\1(h(s))=\pi(\Gamma_x)=\pi(G_x)$ for a point 
$x\in X$, such that $\pi(x)=s.$
\qed \enddemo
\vskip .5cm

\proclaim {Lemma 2.3} In  Case 2  all
 the units in the ring $\CO(S)$ are constants.\ep

\demo{\it Proof of Lemma 2.3} Assume that the Lemma is not true. 
Let $t\in \CO(S)$ be 
a non-constant unit in $\CO(S)$ and $t^*$ be its lift into $\CO(X).$ 
Since there exists a dominant regular map 
$F_I$ (see Lemma1.2) of $\BC^s$ into   any $G$ orbit,  
$t^*$ has to be constant 
along  any $G-$orbit, hence, it has to 
be $G-$ invariant
and provides a  $G-$invariant regular
function $p$,
which is constant over every point $s\in S.$
Hence, it generates a regular function $p_1$ on $S,$
 which is impossible in Case 2.
  \qed \enddemo

\proclaim {Lemma 2.4} In  case 2 $\CO(
S)$ is factorial  
.\ep

\demo{\it Proof of Lemma 2.4}  
 It is enough to show that  any effective
divisor $A$ in $S$ may be defined as
$f=0$ for some $f\in \CO(S).$
Consider once more diagrams (3) and (4), where $\tS$
is any compactification of $S. $   In case 
2 map $p_1$ is not regular and $\tilde {p}_1$ is a 
rational  non-regular map of $\tilde S$ onto $\BP^1.$

At first we are going to prove that $\tilde S$
is rational.

Let $(S',\a)$ a resolution
of $\tS,$ such that the lift $p_1':S' \to \BP^1$ 
of $\tilde p_1$ onto
$S'$ is regular.  Let $B'$ be a normal Stein factorization
of $p_1'.$  We obtain the following
commutative diagram.

$$\alignat {12}\
&\ \ S'  \overset{q}\to\longrightarrow  \ \ B'\\
&\a \downarrow  \  \searrow p_1'\ \ \downarrow  \tau  \tag5 \\
& \ \ \tS  \overset{\tilde p_1}\to\longrightarrow  \ \ \BP^1\endalignat$$ 

In this diagram the maps $p_1', q, \a,\tau$ are regular, $ \tau$
is finite and $\a$ is a blowing down of  finite number 
of exceptional curves (\cite {Sh1},\cite {Sh2}).

 The map $p_1'$ is not constant on the exceptional divisor
 of the map $\a.$ That means that there is an irreducible 
 component $E$ of this divisor, which is mapped by 
 $p_1'$ onto $\BP^1,$ and, hence, $q(E)=B'.$
Since $E$ is exceptional, it has to be rational,
 and $B'$ is rational as well.

Consider now the fiber $F_{b'}=q^{-1}(b')$
over a general point $b'\in B'.$ We have: 
$$ \a(F_{b'})\cap S=\pi (\Gamma_x),$$
where $\Gamma_x$ is an orbit of any point $x\in X,$
such that $\pi(x)\in \a(F_{b'}).$
Since $\Gamma_x \cong \BC$  and $\a|_{F_{b'}}$ is birational,  
it follows that $F_{b'}$ is a rational curve. 

 We obtained that   $S'$ is a ruled surface with rational base.
 Therefore, $S'$ is a rational surface and  $H^1(S',\BZ)=0.$
Since the group $H^1$ is invariant under the blowing-downs, it means
that $H^1(\tS, \BZ)={0}.$
  
   Let $D= \sum D_i=\tS \setminus S.$
Let $A$ be any irreducible curve in $S,$  
and $\tilde A$ be its closure in $\tS.$
Since  $H^2(S,\BZ)={0},$ 
 the map $H_2(D, \BZ) \to H_2(\tS,\BZ)$ is an epimorphism, 
 hence 
 $\tilde A= \sum a_i D_i$ in $H_2 (\tS,\BZ).$ 
 But since $H^1(\tS,\BZ)={0},$
 from the topological equivalence follows the linear equivalence
(see, e.g. \cite {B},p. 7),
 hence   $\tilde A= \sum a_i D_i$ as a divisor, 
 and there exist a function 
 in $\CO(S)$, such that  $A$ is its zero divisor. 
 \qed \enddemo
 \vskip .5 cm
 
 Thus, in  Case 2, $S$ is a smooth affine surface with factorial
 $\CO(S)$ without non-constant units, and $H^2(S,\BZ)=\{0\}.$
 Moreover, the logarithmic Kodaira dimension
 $k(S)=k(\tS\setminus D)=k(S'\setminus \a^{-1}(D))=-\infty, $
 since  the fiber of the restriction of $q$ onto 
 $S'\setminus \a^{-1}(D)$ is isomorphic to $\BC.$ 
By the Miyanishi- Sugie Theorem ([M-S],[Su]),
 $S$ is isomorphic to $\BC^2.$

\vskip .5cm

\proclaim {Lemma 2.5} If $\dim \overline{G_{x}} =n+2$ for a general point 
$x\in X$ , then $S$  is isomorphic to $\BC ^2$ .
\endproclaim

\demo {\it Proof of Lemma 2.5} By virtue of Proposition 1.1 
in this case 
the general orbit contains 
 an image of $\BC^k$ under a regular rational map 
$F_I$ for some $k\geq n+2$ as a dense subset. Let $\overline S,
 \overline X =\overline S\times \BP^n$ be the closures
of $S,X$ respectively, and $ \overline {F}_I: \BP^k 
\to \overline X$ be an extension of $F_I.$ 
Then $\opi\cdot\overline {F}_I$
will be a rational map of $\BP^k$ onto $\overline S.$ 
Since any unirational
surface is rational (\cite {Sh2}. ch 3), $\overline S$ 
has to be rational. That means that,  as in Lemma 2.4, 
from $H^2(S,\BZ)=\{0\} $ follows that $\CO(S)$ is factorial.
Moreover, $\CO(S)$ has no non constant units and $k(S)=-\infty,$ 
since $S$ is dominated by 
$\BC ^k.$ By the Miyanishi -Sugie Theorem ([M-S], [Su]), $S$ is a plane.
\qed\edm  

\demo{\it Proof of Theorem 2.1} 

Let function $f\in \CO(S)$ be invariant under all $\BC-$  
actions on $S.$ We have to prove that 
its  lift $f^*$ onto the product $X$ 
 is invariant under any  
 $\BC-$action  $X$. 
If a general  $G-$orbit is contained  in a fiber of projection $\pi,$
(lies over one point of $S,$) then it is obviously true.
 If $\dim G_x =n+2$ for a general point $x\in X, $ by Lemma 2.5 
 $AK(S)=\BC$
 and the statement of Theorem is valid as well.

Thus, we may assume, that $\dim G_x=n+1$ for an orbit $G_x$ 
of a general point $x\in X.$
Let $f\in AK(S), f\neq const, $ and let  $f^*$ be its
 lift into $\CO(X).$
Let $ x_1,x_2$ be two points in $X$ belonging to the same $G-$ orbit, 
Then $\pi (x_1), \pi(x_2) $, by Lemma 2.2, belong to the same 
orbit of a $\BC-$action on $S,$ thus $f(\pi(x_1)=f(\pi(x_2)).$
 But the lift is invariant under the actions along the fibers of $\pi,
 $ hence $f^*(x_1)=f^*(x_2).$ \qed.\enddemo

\vskip 0.5cm
\proclaim{Remark} In the proof of the Theorem 2.1
 the fiber $\BC^n$ of 
the product $X=S\times \BC^n$
may be replaced by any other affine variety
such that the group of its $\BC-$actions has a
dense orbit.  This is the only property of the fiber used 
in the proof.\ep

\baselineskip 15pt

\noindent Tatiana M. Bandman, Dept. of Mathematics \& CS,
Bar-Ilan University, Ramat-Gan, 52900, Israel,
 e-mail:bandman\@macs.biu.ac.il.

\noindent Leonid Makar-Limanov, Dept. of Mathematics \& CS,
Bar-Ilan University, Ramat-Gan, 52900, Israel, e-mail:
lml\@macs.biu.ac.il; Dept. of Mathematics, Wayne State University,
Detroit, MI 48202, USA, e-mail: lml\@math.wayne.edu. 
\end

\end